\documentclass[12pt,a4paper]{article}
\usepackage{blindtext}
\usepackage{amsmath}
\usepackage{amsfonts}
\usepackage{amsthm}
\usepackage{placeins}
\usepackage{multirow}
\usepackage{xcolor}
\usepackage{float}
\usepackage{graphicx}
\usepackage{caption}
\usepackage{subcaption}
\usepackage{hyperref}
\usepackage[super]{cite}
\usepackage{fullpage}
\usepackage{booktabs}

\newtheorem{theorem}{Theorem}[subsection]

\numberwithin{equation}{section}
\theoremstyle{definition}

\newtheorem{lem}[theorem]{Lemma}

\begin{document}
\begin{center}
\Large \textbf {A hybrid Haar wavelet collocation method for nonlocal hyperbolic partial differential equations }\\
\end{center}
\begin{center}
Gopal Priyadarshi$^{1, 2}$ and Abdul Halim$^{1, 3}$\\
$^1$ Applied Mathematics and Computational Sciences\\
King Abdullah University of Science and Technology (KAUST)\\
Thuwal, 23955--6900, Saudi Arabia\\
\vspace{0.1 cm}
$^2$ Department of Mathematics\\
S.M.D College, Patliputra University, India\\
\vspace{0.1 cm}
$^3$Department of Mathematics\\
H.S. College, Munger University, India
\end{center}
\textbf{Abstract}: In this paper, we propose a hybrid collocation method based on finite difference and Haar wavelets to solve nonlocal hyperbolic partial differential equations. Developing an efficient and accurate numerical method to solve such problem is a difficult task due to the presence of nonlocal boundary condition. The speciality of the proposed method is to handle integral boundary condition efficiently  using the given data.  Due to various attractive properties of Haar wavelets such as closed form expression, compact support and orthonormality, Haar wavelets are efficiently used for spatial discretization and second order finite difference is used for  temporal discretization. Stability and error estimates have been investigated in order to ensure the convergence of the method. Finally, numerical results are compared with few existing results and it is shown that numerical results obtained by the proposed method is better than few existing results. 
\textbf{Mathematics Subject Classification(2010) 65T60, 65C30, 31A30.}\\
\section{Introduction}

Numerical techniques for solving nonlocal hyperbolic partial differential equations have received enormous attention over the last few years. These nonlocal hyperbolic PDEs are used to describe the dynamics of ground water (see \cite{Bouziani1, Pulkina}). Some problems in visco-elasticity and food industry are also described in terms of nonlocal hyperbolic partial differential equations (see \cite{Bouz, Renardy, Kavalloris}). The nonlocal boundary conditions appear in the hyperbolic PDEs when the boundary data cannot be measured directly. 

In this article, we consider the non local hyperbolic problem given by
\label{sec:6.1}
\begin{equation}\label{6.1.1}
\frac{\partial^2 u}{\partial t^2} - \frac{\partial^2 u}{\partial x^2} = \phi(x, t), ~x \in (0, 1), ~t \in [0, T],
\end{equation}
with initial conditions
\begin{equation}\label{6.1.2}
u(x, 0) = f(x), ~~~~~~~~0 \leq x \leq 1,
\end{equation}
\begin{equation}\label{6.1.3}
\frac{\partial u}{\partial t}(x, 0) = g(x), ~~~~~~ 0 \leq x \leq 1,
\end{equation}
Dirichlet boundary condition
\begin{equation}\label{6.1.4}
u(0, t) = h(t), ~~~~~~~~~~ 0 < t \leq T,
\end{equation}
and nonlocal condition
\begin{equation}\label{6.1.5}
\int_{0}^{1}u(x, t)dx = \nu(t), ~~0 < t \leq T,
\end{equation}
where $\phi, f, g, h$ and $\nu$ are known functions. It is assumed that
\[f(x) \in C[0, 1] \cap C^{2}[0, 1] ~\mbox{and}~ g(x) \in C[0, 1] \cap C^{1}[0, 1].\]
Further, we assume that $\phi$ is sufficiently smooth in order to obtain a smooth solution $u$. It should be noted that $f(x)$ and $g(x)$ satisfy the following compatibility conditions
\[f(0) = h(0), \int_{0}^{1}f(x)dx = \nu(0),\]
\[g(0) = h^{'}(0), \int_{0}^{1}g(x)dx = \nu^{'}(0). \]

The existence, uniqueness and stability results for the given problem \eqref{6.1.1}--\eqref{6.1.5} that combine integral as well as Neumann conditions are discussed by Beilin \textit{et al.} \cite{Beilin}. Gordeziani  \textit{et al.} \cite{Gordeziani}  and Kavalloris \textit{et al.}\cite{Kavalloris}  have also investigated hyperbolic partial differential equations with nonlocal boundary conditions.

Numerical techniques for the solution of nonlocal hyperbolic equations have been proposed by several researchers. Ang \textit{et al.} \cite{Ang} proposed a numerical method based on integro-differential equation and local interpolating functions to solve the nonlocal hyperbolic PDEs. Dehghan and his collaborators investigated nonlocal hyperbolic PDEs using several numerical methods, e.g. finite difference method based on cubic B-spline scaling functions \cite{Dehgan7},  variational iterative method \cite{Dehgan12a}, meshless method using radial basis functions \cite{Dehgan8},  and Bernstein Ritz-Galerkin method \cite{Yousefi1}. A numerical method based on shifted Legendre tau technique has been proposed by Saadatmandi \textit{et al.}\cite{Saadatmandi}. 

In the last few years, Haar wavelet based collocation methods are extensively used for the numerical solution of partial differential equations. Because of various attractive properties of Haar wavelet such as closed form expression, compact support and orthonormality, it is widely used in various areas of science and engineering. The main drawback of Haar wavelet is its discontinuity. Therefore, we cannot express the solution in terms of Haar wavelet basis directly. To come out with this difficulty, either we can regularize the Haar wavelets with interpolating splines (see \cite{Cattani}) or expand the highest derivatives in terms of Haar wavelet basis and integrate it out to get the desired expressions (see \cite{Chen}). We have used later approach to  handle the difficulty coming from discontinuity of the wavelet. Second order finite difference scheme is used for the temporal discretization whereas Haar wavelet basis is used for the spatial discretization. Stability and error analysis have been rigorously studied in order to ensure the convergence of the method. The obtained numerical results are compared with the numerical results provided in the paper \cite{Dehgan8}  by Dehghan \textit{et al.} In paper \cite{Dehgan8}, authors reformulated the problem in such a way that  the integral boundary condition is converted into  a periodic boundary condition. We have dealt with integral boundary condition directly using the given data which is  more accurate.

The content of this paper is organized as follows. In section 2, we review some basic background of Haar wavelet. In section 3, we  propose a hybrid Haar wavelet collocation method (HHWCM) for nonlocal hyperbolic PDEs. In section 4, Stability and error analysis have been studied. Numerical results are analyzed in section 5. A brief conclusion is presented in section 6.
\section{Basic Background}
In this section, we review some basics of Haar wavelet which will be used for the proposed numerical method.
\subsection{Haar wavelet}\label{sec:haarwavelet}
For $i \geq 2$, Haar wavelet family $\{h_{i}(x)\}$ is defined as
\begin{equation}\label{Haarmatrix}
 h_{i}(x) = \left\{
  \begin{array}{l l}
    1, & \quad \text{for $\frac{k}{m} \leq x < \frac{k+0.5}{m}$},\\\\
   -1, & \quad \text{for $\frac{k+0.5}{m} \leq x < \frac{k+1}{m}$},\\\\
    0, & \quad \text{elsewhere},
  \end{array} \right.
  \end{equation}
where $m = 2^j, ~j = 0, 1, 2, ..., J$, ~$k = 0, 1, \cdots, m-1$ and $i = m + k + 1$. $i$ and $J$ denote the wavelet number and maximum resolution level respectively.  For simplicity, we have considered $x \in [0, 1]$.\\
It is to be noted that $h_{1}(x)$ correspond to Haar scaling function defined by 
\begin{equation}\label{Haarmatrix}
 h_{1}(x) = \left\{
  \begin{array}{l l}
    1, & \forall x \in [0, 1),\\
    0, & \text{elsewhere}.
  \end{array} \right.
  \end{equation}
Haar wavelets are orthogonal functions i.e.
\begin{equation}\label{9}
\int_{0}^{1}h_{\alpha}(x)h_{\beta}(x)dx = \left\{
  \begin{array}{l l}
    2^{-j}, & \quad \text{if $\alpha = \beta,$}\\
    0, & \quad \text{if $\alpha \neq \beta.$}
  \end{array} \right.
\end{equation}
The wavelet approximation of a function $u \in L^2[0, 1)$ is given by 
\begin{equation}
u_{J}(x) = \sum_{i = 1}^{2M}a_{i}h_{i}(x),
\end{equation}
where $a_{i}$ is the Haar wavelet coefficient.\\
In order to solve PDEs of any order, we need to compute the following integrals
\begin{equation}\label{10}
P_{\beta, i}(x) = \int_{0}^{x} \int_{0}^{x} ... \int_{0}^{x}h_{i}(t)dt^{\beta} = \frac{1}{(\beta -1)!}\ \int_{0}^{x}(x-t)^{\beta -1 }h_{i}(t)dt,
\end{equation}
where $\beta = 1, 2, \cdots, n$ and $i = 1, 2, 3, \cdots, 2M$.
\noindent Using the definition of Haar wavelets, these integrals are calculated as follows:
\begin{equation}\label{Haarintegral1}
P_{\beta, i}(x) = \left\{
  \begin{array}{l l}
    0, & x < \frac{k}{m}, \\
    \frac{1}{\beta!}\Big(x - \frac{k}{m}\Big)^{\beta},  & x \in \Big[\frac{k}{m}, \frac{k+0.5}{m}\Big),\\
    \frac{1}{\beta!}\Big[\Big(x - \frac{k}{m}\Big)^{\beta} - 2\Big(x - \frac{k+0.5}{m}\Big)^{\beta}\Big], & x \in \Big[\frac{k+0.5}{m}, \frac{k+1}{m}\Big),\\
    \frac{1}{\beta!}\Big[\Big(x - \frac{k}{m}\Big)^{\beta} - 2\Big(x - \frac{k+0.5}{m}\Big)^{\beta} + \Big(x - \frac{k+1}{m}\Big)^{\beta}\Big], & x \geq \frac{k+1}{m}.
  \end{array} \right.
\end{equation}
In the special case, when $\beta = 1$ and 2, we obtain
\begin{equation}\label{Haarintegral2}
P_{1, i}(1) = \left\{
  \begin{array}{l l}
    1,   & \quad \text{for $i = 1,$}\\
    0, & \quad \text{for $i \neq 1.$} 
    
    \end{array} \right.
\end{equation}
and
\begin{equation}\label{Haarintegral3}
P_{2, i}(1) = \left\{
  \begin{array}{l l}
   0.5,   & \quad \text{for $i = 1,$}\\
    \displaystyle \frac{1}{4m^2},  & \quad \text{for $i \neq 1.$} 
    \end{array} \right.
\end{equation}
Let us define
\begin{equation}\label{Haarintegral4}
C_{1, i} = \int_{0}^{1} P_{1, i}(x) dx = \left\{
  \begin{array}{l l}
   0.5,   & \quad \text{for $i = 1,$}\\
    \displaystyle \frac{1}{4m^2},  & \quad \text{for $i \neq 1.$} 
    \end{array} \right.
\end{equation}
and 
\begin{equation}\label{Haarintegral5}
C_{2, i} = \int_{0}^{1} P_{2, i}(x) dx = \left\{
  \begin{array}{l l}
   \frac{1}{6},   & \quad \text{for $i = 1,$}\\
    \displaystyle \frac{2m - 2k -1}{8m^3},  & \quad \text{for $i \neq 1.$} 
    \end{array} \right.
\end{equation}
The grid points are given by
\[y_l =  l\Delta y, ~l = 0, 1, 2, \cdots, 2M.\]
where $\Delta y = \frac{1}{2M}$.\\
The collocation points are given as 
\[x_l = \frac{y_l + y_{l-1}}{2}, ~l = 1, 2, \cdots, 2M.\] 
\noindent Next, we introduce Haar matrix, $H,$ and Haar integral matrices $P_{1}$ and $P_{2}$ which are square matrices of size $2M \times 2M.$ The elements of these matrices are $H(i, l) = h_{i}(x_l), P_{1}(i, l) = P_{1, i}(x_l)$ and  $P_{2}(i, l) = P_{2, i}(x_l)$\\\\
\textbf{Temporal discretization}:
Let $T_{final}$ be the final time where we want to compute the solution. The temporal discretization is given by:
\[0 = t_0 < t_1 < t_2, \dots, t_r, \dots <t_N = T_{final}\]
where $r = 0, 1, 2, \dots, N$ and $t_r = rdt$ where $dt = \frac{T_{final}}{N}$.
\section{A hybrid Haar wavelet collocation method for nonlocal hyperbolic partial differential equation}
In this section, we propose a hybrid wavelet collocation method based on Haar wavelets and second order finite difference method to solve the problem \eqref{6.1.1}--\eqref{6.1.5}. We assume that $u_J$ be the wavelet approximation of $u.$\\\\
Let us assume
\begin{align}\label{haar}
   \frac{\partial^2 u_J}{\partial x^2} (x,t)=\sum\limits_{i=1}^{2M} a_i(t) h_i(x)
\end{align}
Integrating equation \eqref{haar} from $0$ to $x$, we get
\begin{align}\label{haar1}
    \frac{\partial u_{J}}{\partial x} (x,t)=\sum\limits_{i=1}^{2M} a_i(t) P_{1,i}(x) + \frac{\partial u_{J}}{\partial x} (0,t)
\end{align}
Integrating equation \eqref{haar1} from $0$ to $1$, we get
\begin{align}
    u_{J}(1,t)-u_{J}(0,t)= \sum\limits_{i=1}^{2M} a_i(t) C_{1,i} + \frac{\partial u_{J}}{\partial x} (0,t)
\end{align}
where \[C_{1, i} = \int_{0}^{1} P_{1, i}(x) dx.\]
Hence,
\begin{align}\label{haar2}
    \frac{\partial u_J}{\partial x} (x,t)=\sum\limits_{i=1}^{2M} a_i(t) \Big(P_{1,i}(x) - C_{1, i}\Big) + u_J(1,t)-u_{J}(0,t)
\end{align}

\noindent Integrating equation (\ref{haar2}) from $0$ to $x$, we get
\begin{align}\label{haar3}
    u_J(x,t)= \sum\limits_{i=1}^{2M}  a_i(t) \Big (P_{2,i}(x)-xC_{1,i}\Big) + x [u_J(1,t) - u_J(0,t)] + u_J(0,t)
\end{align}
Using nonlocal condition (\ref{6.1.5}), we obtain

\begin{align}
    \sum\limits_{i=1}^{2M}  a_i(t) \Big(C_{2,i}-\frac{1}{2}C_{1,i}\Big) + \frac{1}{2} [u_J(1,t)-u_J(0,t)] + u_J(0, t) = \nu(t)
\end{align}
After  simplification, we get
\begin{align}\label{haar4}
   u_J(1,t)= \sum\limits_{i=1}^{2M} a_i(t) \Big( C_{1, i}- 2C_{2,i}\Big) + 2\nu(t) - h(t)
\end{align}
Thus, from equation \eqref{haar3}, we get
\begin{align}\label{haar51}
    u_J(x,t)&=\sum\limits_{i=1}^{2M} a_i(t)\Big(P_{2, i}(x) -2xC_{2, i}\Big) + 2x[\nu(t) - h(t)] + h(t)
\end{align}
Using second order finite difference scheme for temporal discretization and Haar wavelets for spatial discretization, we obtain
\begin{align}\label{dis_eqn}
\frac{u_{J}(x, t_{n+1}) - 2u_{J}(x, t_{n}) + u_{J}(x, t_{n-1})}{\Delta t^2} = \sum\limits_{i=1}^{2M} a_i h_i(x) + \phi (x, t_{n})
\end{align}
From the given boundary condition (\ref{6.1.3}) and using central difference formula, we get
\begin{align}
\frac{u_{J}(x, t_1) - u_{J}(x, t_{-1})}{2 \Delta t} = g(x)
\end{align}
This implies
\begin{align}\label{ic_dis}
u_J(x, t_{-1}) =  u_J(x, t_1) - 2 \Delta t  g(x)
\end{align}
Using \eqref{dis_eqn} and \eqref{ic_dis}, we obtain the following equation at $t_{0} = 0$
\begin{align}\label{haar61}
u_J(x, t_1) = u_J(x, t_0) + \Delta t g(x) + \frac{(\Delta t)^2}{2} \sum\limits_{i=1}^{2M} a_i h_i(x) +  \frac{(\Delta t)^2}{2}\phi (x, t_0)
    \end{align}
Using equation (\ref{haar51}) and (\ref{haar61}), we obtain
    \begin{align}
    \nonumber
    &\sum\limits_{i=1}^{2M} a_i\Big(P_{2, i}(x) - \frac{\Delta t^2}{2}h_{i}(x) - 2xC_{2, i} \Big) + 2x[\nu(t_1) - h(t_1)] + h(t_1) \\&= u_J(x, t_0) + \Delta t g(x)  +  \frac{\Delta t^2}{2}\phi (x, t_0)
\end{align}
At $t = t_n$, we obtain the following discretized scheme,
\begin{align}\label{final_ut}
\nonumber
    &\sum\limits_{i=1}^{2M} a_i \Big(P_{2, i}(x) - \Delta t^2 h_{i}(x) - 2xC_{2, i} \Big) + 2x[\nu(t_{n+1}) - h(t_{n+1})] + h(t_{n+1}) \\&=  \Delta t^2 \phi(x, t_n) + 2u_J(x, t_n) - u_J(x, t_{n-1})
\end{align}
Equation (\ref{final_ut}) at the collocation points $x_{l}, l = 1, 2, 3, \dots, 2M$ is given by
\begin{align}\label{final_formula}
\nonumber
    &\sum\limits_{i=1}^{2M} a_i \Big(P_{2, i}(x_l) - \Delta t^2 h_{i}(x_l) - 2x_lC_{2, i} \Big) + 2x_l[\nu(t_{n+1}) - h(t_{n+1})] + h(t_{n+1}) \\&=  \Delta t^2 \phi(x_l, t_n) + 2u_J(x_l, t_n) - u_J(x_l, t_{n-1})
\end{align}
Finally, we obtain a matrix system at $t = t_{n}$
\[\textbf{Ba = c}\]
where $\textbf{B} = \{b_{li}, ~1 \leq l, i \leq 2M\}$ and $\textbf{c} = \{c_{l}, ~1 \leq l \leq  2M\}$. The expression for $b_{li}$ and $c_l$ is given by 
\[b_{li} = \Big(P_{2, i}(x_l) - \Delta t^2 h_{i}(x_l) - 2x_lC_{2, i} \Big)\] 
and
\[c_l = \Delta t^2 \phi(x_l, t_n) + 2u_J(x_l, t_n) - u_{J}(x_l, t_{n-1}) -2x_l[\nu(t_{n+1}) - h(t_{n+1})] - h(t_{n+1})\]
At each time step, we calculate the wavelet coefficients $\textbf{a}$ and obtain the required solution.






\section{Stability and error analysis}
\subsection{Stability analysis}
In this subsection, we will study stability analysis for the hybrid Haar wavelet collocation method.\\
Equation \eqref{6.1.1} can be written as follows:
\begin{equation}
\frac{\partial^2 u}{\partial t^2}(x, t) =  \mathcal {L} u(x, t) + \phi(x, t)
\end{equation}
where $\mathcal {L} = \frac{\partial^2 }{\partial x^2}$ is the differential operator. Following the temporal discretization using finite difference technique, we obtain
\begin{equation}\label{disc_eq}
u^{n+1} - 2u^{n} + u^{n-1} = (\Delta t) ^2 \mathcal {H} u^{n+1} + (\Delta t) ^2 \phi(x, t_{n})
\end{equation}
\begin{equation}
\implies u^{n+1} = 2(I - (\Delta t) ^2 \mathcal {H})^{-1} u^{n}  - (I - (\Delta t) ^2 \mathcal {H})^{-1} u^{n-1}  + (I - (\Delta t) ^2 \mathcal {H})^{-1} (\Delta t) ^2 \phi(x, t_{n})
\end{equation}
where $I$ is the identity matrix and $\mathcal {H}$ is the Haar matrix corresponding to the differential operator $\mathcal {L}$.

Since equation involves two time levels, we add one identity equation in order to make single time level. We proceed as follows:
\begin{align}\label{haar71}
u^{n+1} &= B_1 u^{n} + B_2 u^{n-1}\\\nonumber
u^{n} &= u^{n}
    \end{align}
    where $B_1 = 2(I - (\Delta t) ^2 \mathcal {H})^{-1}$ and $B_2 = -(I - (\Delta t) ^2 \mathcal {H})^{-1}.$ Equation (\ref{haar71}) can be written in the  matrix form as follows
    \begin{equation}
    \begin{bmatrix}
B_1 & B_2\\
I_{2M} & \pmb{0}_{2M}
\end{bmatrix}
 \begin{bmatrix}
u^{n}\\
u^{n-1}
\end{bmatrix}
= \begin{bmatrix}
u^{n+1}\\
u^{n}
\end{bmatrix}
 \end{equation}
 We know that the eigenvalues of identity matrix is always 1. The stability of the numerical scheme will depend upon the eigenvalues of the matrix $B$ where 
\begin{equation}
     B =
 \begin{bmatrix}
B_1 & B_2\\
I_{2M} & \pmb{0}_{2M}
\end{bmatrix}
\end{equation}
The proposed numerical scheme will be stable if all the eigenvalues of the matrix $B$ is less than or equal to 1.\\\\
The eigenvalues of the matrix $B$ for hybrid Haar wavelet collocation method for different $\Delta t$ and $J$ are given below
\begin{figure}[H]
\includegraphics[width = 0.5\textwidth, height = 0.4\textwidth]{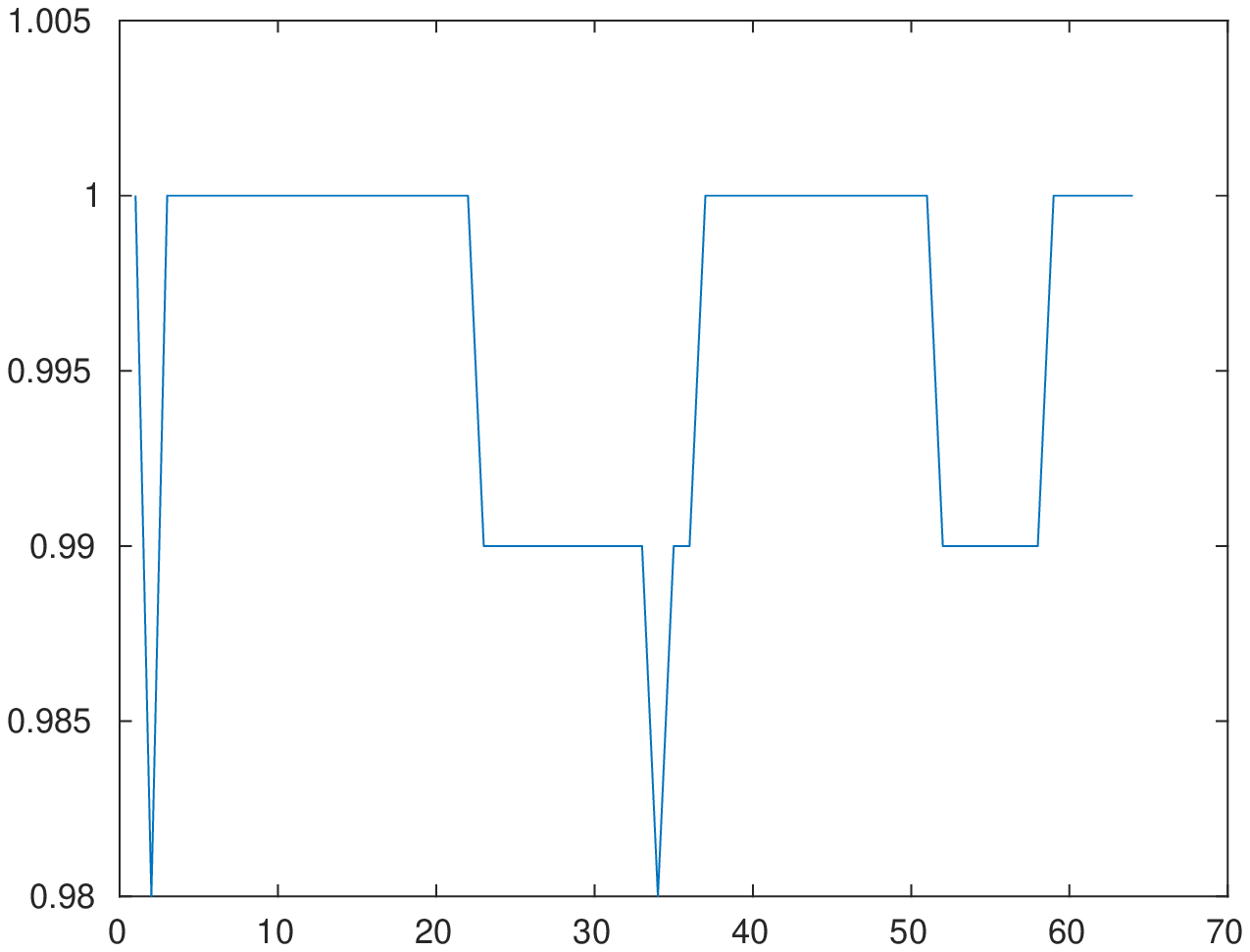}
\hspace{0.4 cm}
 \includegraphics[width = 0.5\textwidth, height = 0.4\textwidth]{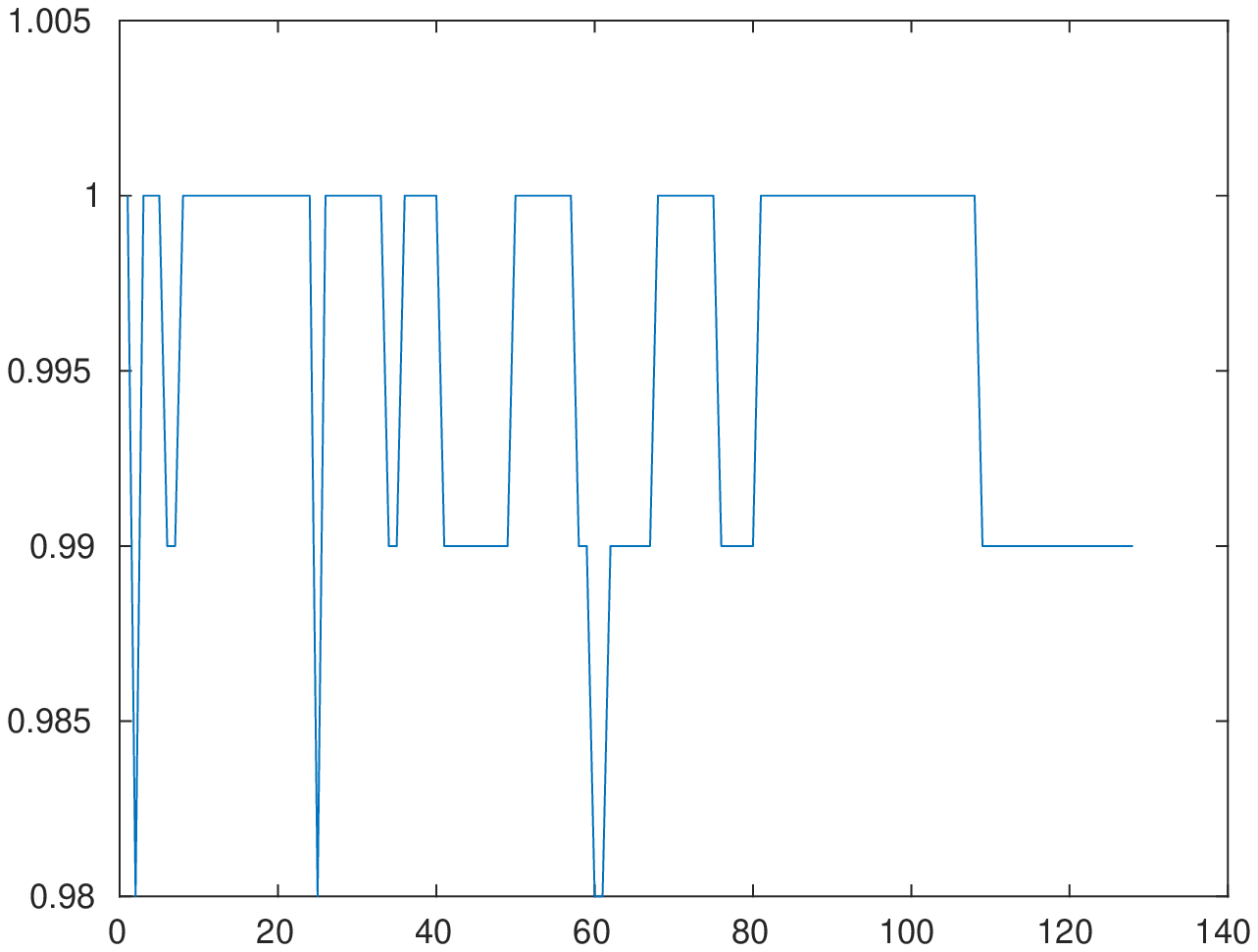}
\caption{Eigenvalues of $B$ at $\Delta t = 10^{-2}$ and $J = 4, 5.$ }
\label{fig3}
\end{figure}
\begin{figure}[H]
\includegraphics[width = 0.5\textwidth, height = 0.4\textwidth]{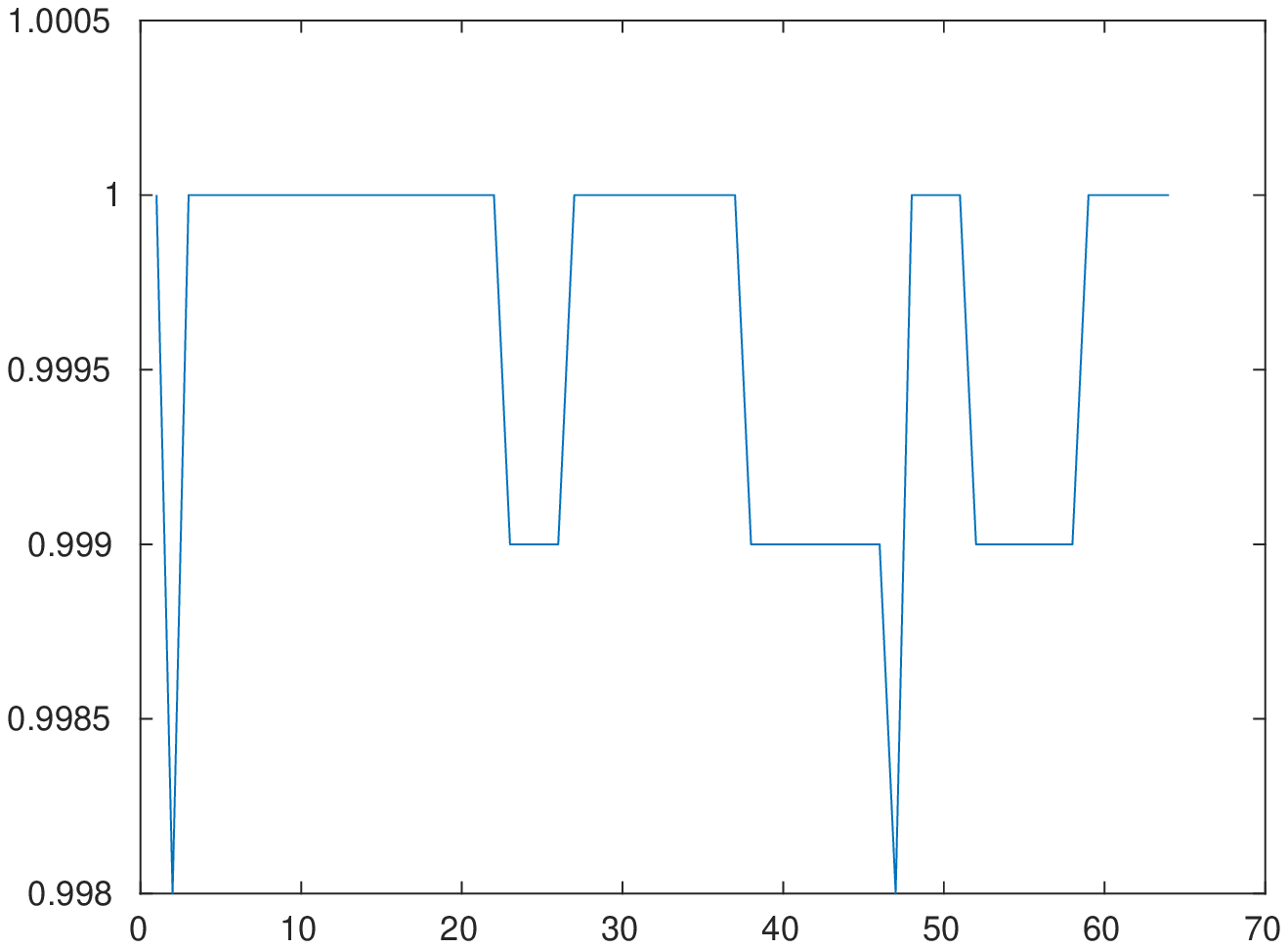}
\hspace{0.4 cm}
 \includegraphics[width = 0.5\textwidth, height = 0.4\textwidth]{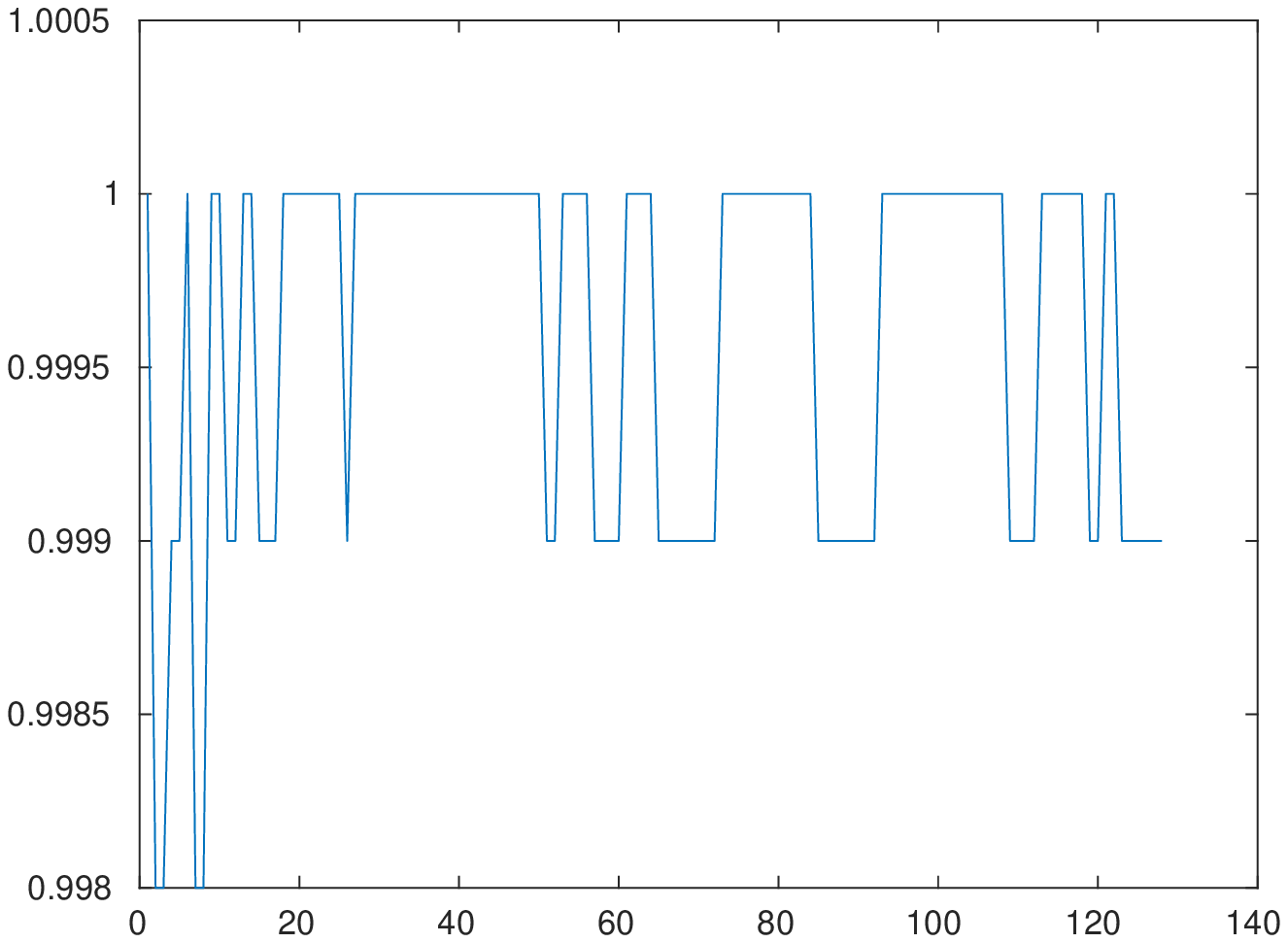}
\caption{Eigenvalues of $B$ at $\Delta t = 10^{-3}$ and $J = 4, 5.$ }
\label{fig3}
\end{figure}
From the above figures, it can be easily guaranteed that the proposed method is stable.
\subsection{Error Analysis}
In this subsection, we will study  error analysis for the proposed numerical method.
From equation \eqref{haar51}, we have the approximate representation of function $u$  given by
\begin{align}\label{haar61a}
    u_{J}(x,t)&=\sum\limits_{i=1}^{2M} a_i\Big(P_{2, i}(x) -2xC_{2, i}\Big) + 2x(\nu(t) - h(t)) + h(t)
\end{align}
and the exact representation of $u$ is given by
\begin{align}\label{haar5}
    u(x,t)&=\sum\limits_{i= 1}^{\infty} a_i\Big(P_{2, i}(x) -2xC_{2, i}\Big) + 2x(\nu(t) - h(t)) + h(t)
\end{align}
Hence, the error term is given by
\begin{align}\label{haar5}
     E_J = u - u_J =\sum\limits_{i= 2M + 1}^{\infty} a_i\Big(P_{2, i}(x) -2xC_{2, i}\Big) 
\end{align}
In terms of resolution level, equation \eqref{haar5} can be written as
\begin{align}\label{haar696}
     E_J = \sum\limits_{j= J + 1}^{\infty} \sum\limits_{k= 0}^{2^{j} - 1} a_{2^{j} + k + 1}\Big(P_{2, 2^{j} + k + 1}(x) -2xC_{2, 2^{j} + k + 1 }\Big) 
\end{align}
\begin{lem}{\cite{Gopal}}
Let us assume that $u$ is a Lipschitz continuous function in [0, 1]. Then, the wavelet coefficient satisfy the following:
\begin{equation}
|a_{2^{j} + k + 1}| \leq \frac{L}{2^{j + 1}}
\end{equation}
where $L \geq 0$ is the Lipschitz constant and 
\begin{equation}
a_{2^{j} + k + 1} = 2^{j} \int_{0}^{1} u(x) h_{2^{j} + k + 1 } (x)dx
\end{equation}
\end{lem}
\begin{lem}
Let $u$ be the Lipschitz continuous function in the unit square. Then, for fixed $t$, the proposed method is convergent and order of convergence is 2 in spatial variable i.e. 
\begin{equation}
||E_{J}||_{2} = ||u - u_{J}||_{2} = \mathcal{O}\Big[\Big(\frac{1}{2^{J} + 1 }\Big)^{2}\Big]
\end{equation}
\end{lem}

\noindent \textbf{Proof} From equation (\ref{haar696}) and definition of $L^2$ norm, we have
\begin{align}\label{7.2}
\Vert E_{J} \Vert^2 & = \int_{0}^{1}   \bigg|\sum_{{j} = J +1}^{\infty}\sum_{k = 0}^{2^{j} - 1} a_{2^{j} + k + 1} \Big (P_{2, 2^{j} + k + 1}(x) - 2xC_{2, 2^{j} + k + 1}\Big)\bigg|^2dx  \nonumber \\
& = \int_{0}^{1} \bigg|\sum_{{j} = J +1}^{\infty}\sum_{k = 0}^{2^{j} - 1}\sum_{{j_1} = J +1}^{\infty}\sum_{k_1 = 0}^{2^{j_1} - 1}a_{2^{j} + k + 1}a_{2^{j_1} + k_1 + 1}
\Big(P_{2, 2^{j} + k + 1}(x)  - 2xC_{2, 2^{j} + k + 1}\Big) \times \nonumber \\ & \hspace{5 cm}\Big(P_{2, 2^{j_1} + k_1 + 1}(x) -2x C_{2, 2^{j_1} + k_1 + 1}\Big)\bigg|dx \nonumber \\
& \leq  \sum_{{j} = J +1}^{\infty}\sum_{k = 0}^{2^{j} - 1}\sum_{{j_1} = J +1}^{\infty}\sum_{k_1 = 0}^{2^{j_1} - 1} |a_{2^{j} + k + 1}| |a_{2^{j_1} + k_1 + 1}|\bigg |\int_{0}^{1}\Big (P_{2, 2^{j} + k + 1}(x)  - 2xC_{2, 2^{j} + k + 1}\Big) \times \nonumber \\ & \hspace{5 cm} \Big (P_{2, 2^{j_{1}} + k_{1} + 1}(x)  - 2xC_{2, 2^{j_{1}} + k_{1} + 1}\Big )\bigg|dx \nonumber \\
& \leq  \sum_{{j} = J +1}^{\infty}\sum_{k = 0}^{2^{j} - 1}\sum_{{j_1} = J +1}^{\infty}\sum_{k_1 = 0}^{2^{j_1} - 1} |a_{2^{j} + k + 1}| |a_{2^{j_1} + k_1 + 1}|\bigg |\int_{0}^{1}\Big (P_{2, 2^{j} + k + 1}(x)P_{2, 2^{j_{1}} + k_{1} + 1}(x)  \nonumber \\ & - 2xC_{2, 2^{j} + k + 1} P_{2, 2^{j_{1}} + k_{1} + 1}(x) - 2xC_{2, 2^{j_{1}} + k_{1} + 1}P_{2, 2^{j} + k + 1}(x) + 4x^2C_{2, 2^{j + k + 1}}C_{2, 2^{j_{1}} + k_{1} + 1}\Big) dx\bigg|
\end{align}
By the above lemma
\begin{equation}\label{7.3}
|a_{2^{j} + k + 1}| \leq\frac{A_1}{{2^{j + 1}}},
\end{equation}
and 
\begin{equation}\label{7.4}
|a_{2^{j_1} + k_1 + 1}| \leq\frac{A_2}{{2^{j_1 + 1}}}.
\end{equation}
\begin{align*}
\bigg|\int_{0}^{1}P_{2, 2^{j} + k + 1}(x)P_{2, 2^{j_1} + k_1 + 1}(x)dx\bigg| &\leq \Vert P_{2, 2^{j} + k + 1} \Vert _{L^2}\Vert P_{2, 2^{j_1} + k_1 + 1} \Vert _{L^2}\nonumber \\& 
\end{align*}
We know that $P_{2, 2^{j} + k + 1}(x)$ is monotonically increasing function in [0, 1], hence
\[ P_{2, 2^{j} + k + 1}(x) \leq P_{2, 2^{j} + k + 1}(1) = \bigg(\frac{1}{{2^{j+1} }}\bigg)^2\]
Using the fact that our domain is of  finite measure, we have the following estimate
\begin{equation*}\label{7.5}
\Vert P_{2, 2^{j} + k + 1}\Vert_{L^2} \leq \bigg(\frac{1}{{2^{j+1} }}\bigg)^2
\end{equation*}
Hence,
\begin{align}\label{7.5a}
\bigg|\int_{0}^{1}P_{2, 2^{j} + k + 1}(x)P_{2, 2^{j_1} + k_1 + 1}(x)dx\bigg|\leq \bigg(\frac{1}{{2^{j+1} }}\bigg)^2 \bigg(\frac{1}{{2^{j_1+1} }}\bigg)^2
\end{align}
Similarly
\begin{align}\label{7.5b}
\bigg|\int_{0}^{1}2xC_{2, 2^{j} + k + 1}P_{2, 2^{j_1} + k_1 + 1}(x)dx\bigg|\leq \bigg(\frac{1}{{2^{j+1} }}\bigg)^2 \bigg(\frac{1}{{2^{j_1+1} }}\bigg)^2
\end{align}
\begin{align}\label{7.5c}
\bigg|\int_{0}^{1}2xC_{2, 2^{j_1} + k_1 + 1}P_{2, 2^{j} + k + 1}(x)dx\bigg|\leq \bigg(\frac{1}{{2^{j+1} }}\bigg)^2 \bigg(\frac{1}{{2^{j_1+1} }}\bigg)^2
\end{align}
and
\begin{align}\label{7.5d}
\bigg|\int_{0}^{1}4x^2C_{2, 2^{j} + k + 1}C_{2, 2^{j_1} + k_1 + 1}dx\bigg|\leq \frac{4}{3}\bigg(\frac{1}{{2^{j+1} }}\bigg)^2 \bigg(\frac{1}{{2^{j_1+1} }}\bigg)^2
\end{align}
Substituting the above estimates (\ref{7.3}--\ref{7.5d}) in (\ref{7.2}), we get
\begin{align*}
\Vert E_{J} \Vert^2 &\leq  \frac{13}{3} A_1A_2 \sum_{{j} = J +1}^{\infty} \sum_{{j_1} = J +1}^{\infty} \bigg(\frac{1}{{2^{j + 1}}}\bigg)^3\bigg(\frac{1}{{2^{j_1 + 1}}}\bigg)^3 2^{j}2^{j_1} \\ & \leq A \sum_{{j} = J +1}^{\infty} \sum_{{j_1} = J +1}^{\infty} \bigg(\frac{1}{{2^{j + 1}}}\bigg)^{2}\bigg(\frac{1}{{2^{j_1 + 1}}}\bigg)^{2},  ~~\mbox{where} \quad (A = \frac{13}{12} A_1A_2)\\&\leq A\bigg [\sum_{{j} = J +1}^{\infty}  \bigg(\frac{1}{{2^{j + 1}}}\bigg) \bigg]^4 \\&  \leq A\frac{1}{2^{4}}\bigg [\frac{1}{2^{J+1}}\sum_{{j} = 0}^{\infty}  \bigg(\frac{1}{{2^{j}}}\bigg) \bigg]^4 \\&
\leq A\bigg(\frac{1}{2^{J+1}}\bigg)^4 ~~~~(\mbox{since} \sum_{{j} = 0}^{\infty}  \bigg(\frac{1}{{2^{j}}}\bigg) = 2.)
\end{align*}
Hence 
\begin{align*}
\Vert E_{J} \Vert\leq K\bigg(\frac{1}{{2^{J + 1}}}\bigg)^2
\end{align*}
where $\displaystyle K = \sqrt{A}$.\\
Hence the order of convergence of the Haar wavelet method in spatial variable is given by
\begin{equation*}
\Vert \mbox{E}_J \Vert _{L_{2}}= \mathcal{O}\bigg[\bigg(\frac{1}{2^{J+1}}\bigg)^2\bigg] 
\end{equation*}
\begin{theorem}
Let $\frac{\partial^2 u}{\partial t^2}$ and $\frac{\partial^2 u}{\partial x^2}$ exist and bounded in $[0, 1] \times [0, T].$ Then the error estimate for the fully discretized hybrid Haar wavelet collocation method is given by  
\begin{equation*}
\Vert \mbox{E} \Vert _{L_{2}}= \mathcal{O}\bigg[\bigg(\frac{1}{2^{J+1}}\bigg)^2 + {\Delta t^2}\bigg].
\end{equation*}
\end{theorem}
\noindent Proof. From the above lemma
\begin{equation*}
||E_{J}||_{2}  = \mathcal{O}\Big[\Big(\frac{1}{2^{J} + 1 }\Big)^{2}\Big]
\end{equation*}
As we have used second order finite difference method for the temporal discretization, the error estimate for the fully discretized numerical method is given by 
\begin{equation*}
\Vert E \Vert _{L_{2}}= \mathcal{O}\bigg[\bigg(\frac{1}{2^{J+1}}\bigg)^2 + {\Delta t^2}\bigg].
\end{equation*}

\section{Results of numerical experiments}
Following the hybrid  Haar wavelet collocation method proposed in section 3, we solve the problem \eqref{6.1.1}-- \eqref{6.1.5} on MATLAB. We present various numerical examples and compare it with few existing results. Our numerical results are better than the existing results \cite{Dehgan8}. \\
\noindent \textbf{Example 1.}
\begin{equation}\label{6.1.20}
\frac{\partial^2 u}{\partial t^2} - \frac{\partial^2 u}{\partial x^2} = \Big(\frac{1}{4} + \pi^2\Big)e^{-\frac{t}{2}}\sin(\pi x), ~~~~0 < x < 1, 0 < t \leq T,
\end{equation}
with initial conditions
\begin{equation}
u(x, 0) = \sin(\pi x), ~~~~~~~~~~ 0 \leq x \leq 1,
\end{equation}
\begin{equation}
\frac{\partial u}{\partial t}(x, 0) = -\frac{1}{2}\sin(\pi x),~~~ 0 \leq x \leq 1,
\end{equation}
and Dirichlet boundary condition
\begin{equation}
u(0, t) = 0, ~~~~~~~~~~~~~~~~~~~~~ 0 < t \leq T,
\end{equation}
with nonlocal condition
\begin{equation}\label{6.1.21}
\int_{0}^{1}u(x, t)dx = \frac{2}{\pi}e^{-\frac{t}{2}}, ~~~~~~~~ 0 < t \leq T.
\end{equation}
The exact solution of (\ref{6.1.20} -- \ref{6.1.21}) is 
\begin{equation*}
u(x, t) = e^{-\frac{t}{2}}\sin(\pi x).
\end{equation*}
Figure \ref{fig3} presents the exact and approximate solutions by the proposed method at different spatial and temporal points. Point wise absolute error at time $T = 1$ and max norm error at different time steps are reported in Table \ref{Table:1} and Table \ref{Table:2} respectively. 
\begin{figure}[H]
\includegraphics[width = 0.5\textwidth, height = 0.4\textwidth]{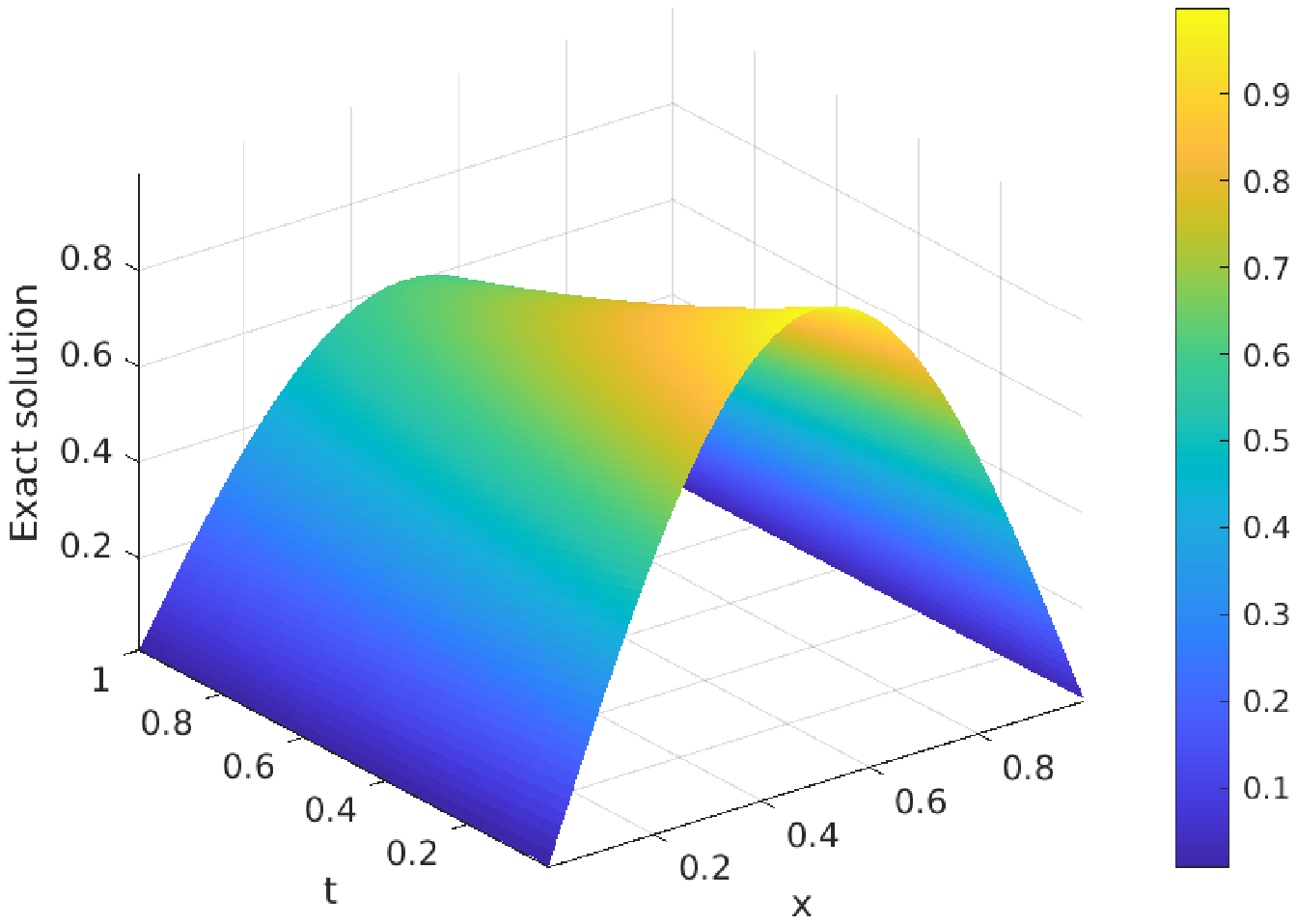}
\hspace{0.4 cm}
 \includegraphics[width = 0.5\textwidth, height = 0.4\textwidth]{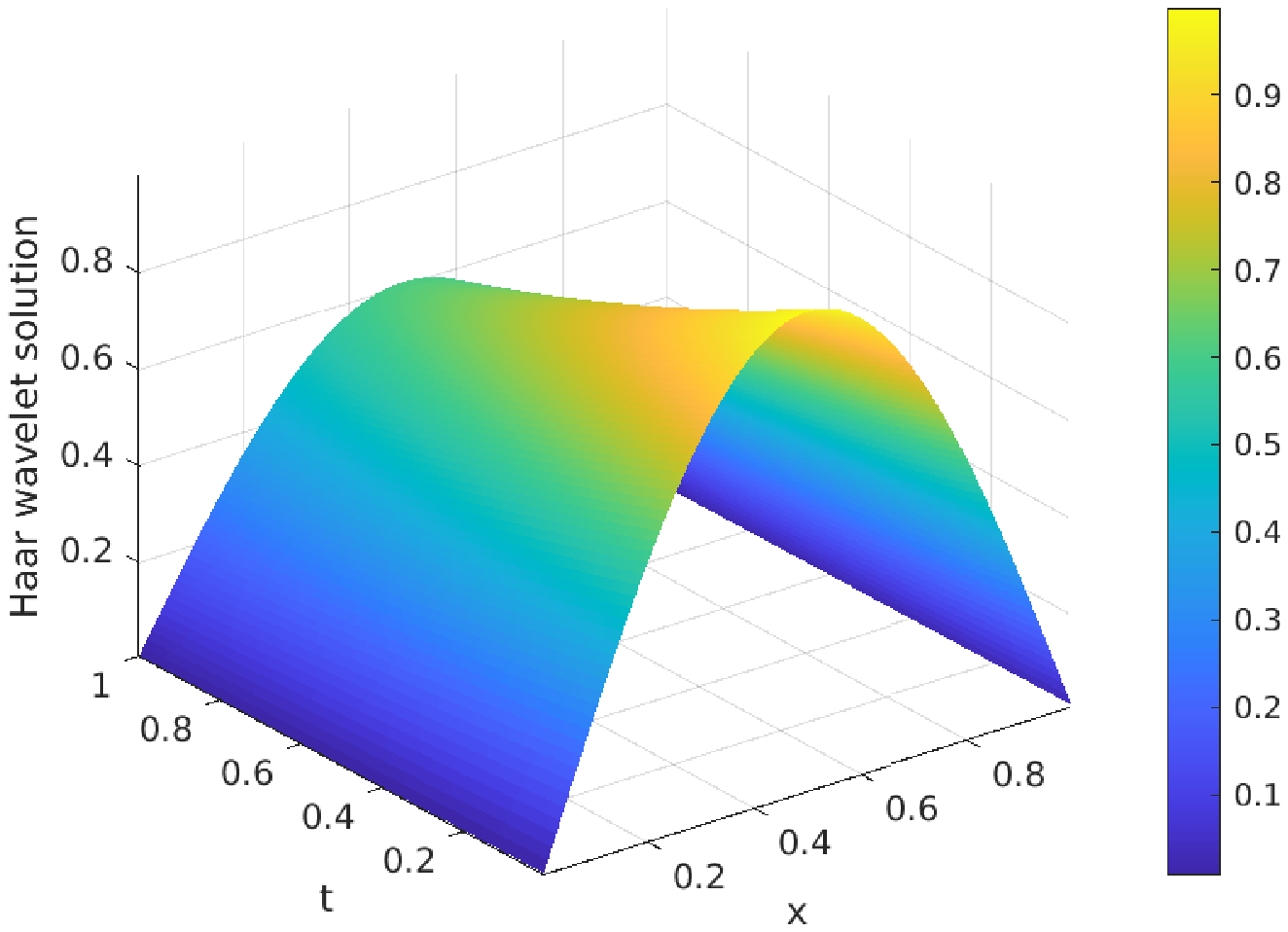}
\caption{Exact and HHWCM based approximate solutions at $J = 6, T = 1$ and $\Delta t = 10^{-4}.$ }
\label{fig3}
\end{figure}
\begin{table}[H]
\centering

\begin{tabular}{|c|c|c|c|c|c|} 
\hline
$x$ & Exact $u$ & Absolute error \\
\hline

0.1  & 0.18742828 & $1.6 \times 10^{-6}$\\ 
0.2  & 0.35650978 & $2.6 \times 10^{-6}$  \\
0.3 & 0.49069361 & $3.0 \times 10^{-6}$   \\ 
0.4   & 0.57684494 & $3.0 \times 10^{-6}$ \\       
0.5  & 0.60653066 & $2.5 \times 10^{-6}$  \\
0.6  & 0.57684494 & $1.5 \times 10^{-6}$ \\
0.7 & 0.49069361 & $3.0 \times 10^{-7}$  \\
0.8 &0.35650978 & $2.9 \times 10^{-6}$ \\
0.9 & 0.18742828& $6.4 \times 10^{-6}$ \\
1.0 & 0.00000000& $8.7 \times 10^{-6}$ \\ 
\hline

\end{tabular}
\caption{Pointwise absolute error at $T = 1, \Delta t = 10^{-4}$ and $J$ = 6. }
\label{Table:1}
\end{table}
\begin{table}[H]
\centering

\begin{tabular}{|c|c|c|c|c|} 
\hline
$T$ & HHWCM &  TPS-RBF \cite{Dehgan8} & MQ-RBF \cite{Dehgan8} & CS-RBF \cite{Dehgan8} \\
\hline
$0.5$ & $3.4 \times 10^{-5}$ & $3.8 \times 10^{-3}$ & $1.3 \times 10^{-3}$  & $2.8 \times 10^{-2}$\\
\hline
$1.0$  & $1.0 \times 10^{-5}$ & $6.8 \times 10^{-3}$ & $2.4 \times 10^{-3}$ & $5.1 \times 10^{-2}$\\ 

\hline
\end{tabular}
\caption{Comparison of maximum error using various numerical methods. }
\label{Table:2}
\end{table}
From the above results, it is observed that a very good accuracy can be achieved at very less resolution level. It is also noticed that maximum absolute error decreases  significantly with reducing $\Delta t$ size. Comparison table shows that the proposed method is better than various meshless method developed by Dehghan \textit{et al.}\cite{Dehgan8} in terms of maximum error.\\
\noindent \textbf{Example 2.}

\begin{equation}\label{6.1.22}
\frac{\partial^2 u}{\partial t^2} - \frac{\partial^2 u}{\partial x^2} = 0, ~~~~0 < x < 1, 0 < t \leq T,
\end{equation}
with initial conditions
\begin{equation}
u(x, 0) = \cos(\pi x), ~~~~~~~~~~ 0 \leq x \leq 1,
\end{equation}
\begin{equation}
\frac{\partial u}{\partial t}(x, 0) = 0,~~~~~~~~~~~~~~~~ 0 \leq x \leq 1,
\end{equation}
and Dirichlet boundary condition
\begin{equation}
u(0, t) =  \cos(\pi t), ~~~~~~~~~~~~ 0 < t \leq T,
\end{equation}
with nonlocal condition
\begin{equation}\label{6.1.23}
\int_{0}^{1}u(x, t)dx = 0, ~~~~~~~~~~~~ 0 < t \leq T.
\end{equation}
The exact solution of (\ref{6.1.22} -- \ref{6.1.23}) is 
\begin{equation*}
u(x, t) = \frac{1}{2}[\cos \pi (x +t) + \cos \pi (x - t)].
\end{equation*}

\begin{figure}[H]
\includegraphics[width = 0.5\textwidth, height = 0.4\textwidth]{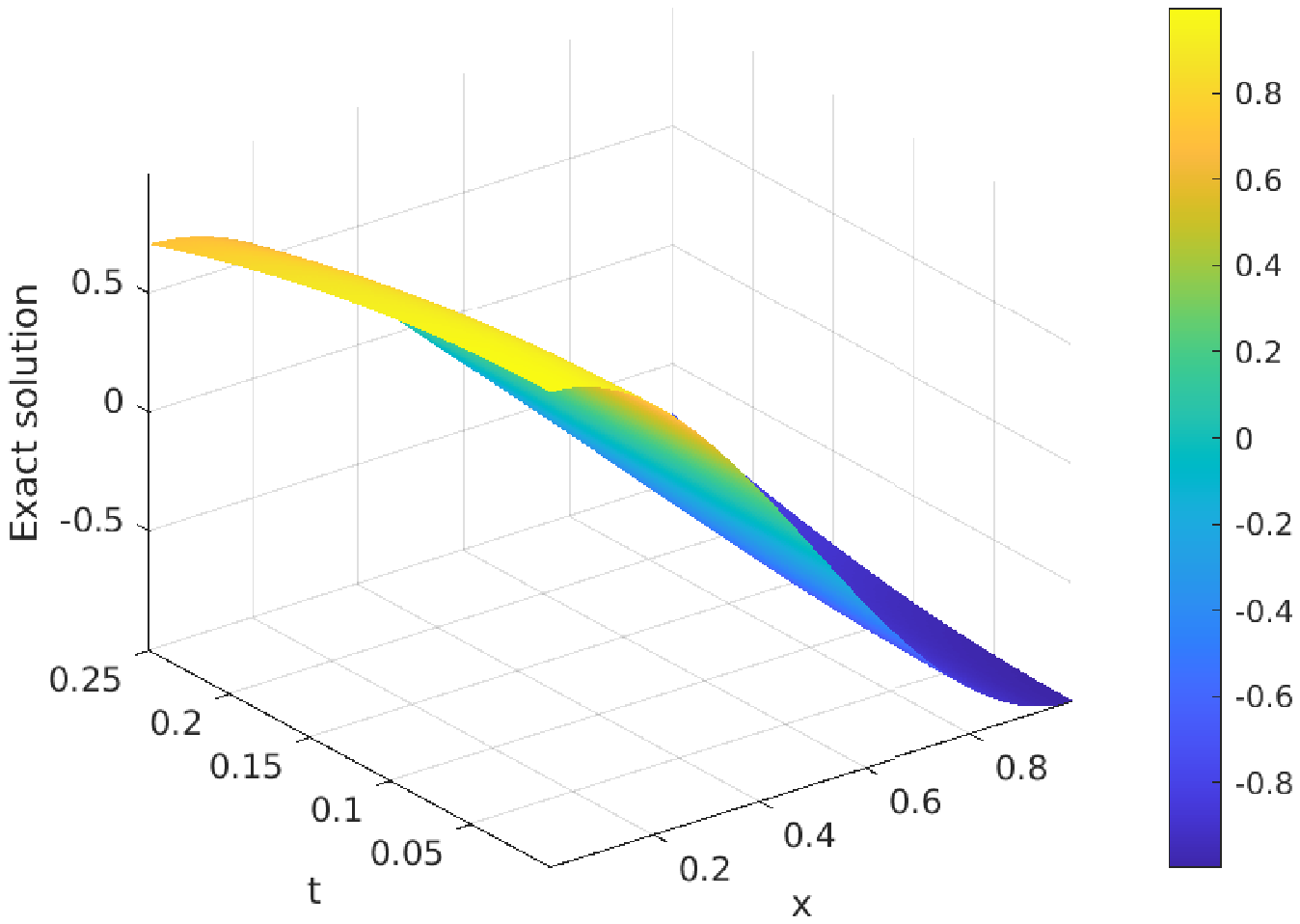}
\hspace{0.4 cm}
 \includegraphics[width = 0.5\textwidth, height = 0.4\textwidth]{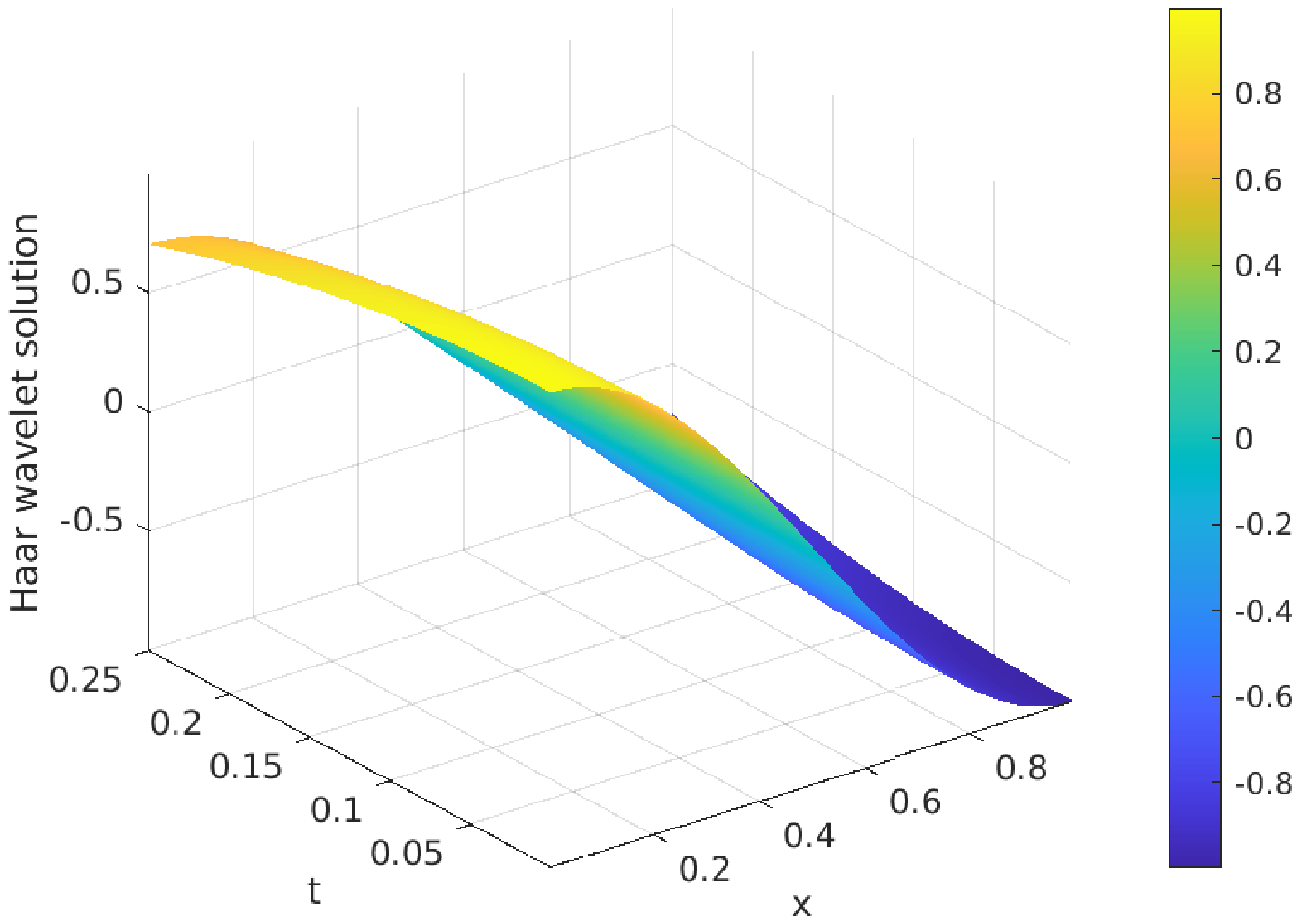}
\caption{Exact and HHWCM based approximate solutions at $J = 6, T = 0.25$ and $\Delta t = 10^{-4}.$ }
\label{fig1}
\end{figure}

\noindent Figure \ref{fig1} presents the exact and approximate solutions by the proposed method at different spatial and temporal points. Point wise absolute error at time $T = 0.25$ and maximum absolute error at different time steps are reported in Table \ref{Table:3} and Table \ref{Table:4} respectively. 
\begin{table}[H]
\centering

\begin{tabular}{|c|c|c|c|c|c|c|c|} 
\hline
$x$ & Exact $u$ & HHWCM  & TPS-RBF \cite{Dehgan8} & MQ-RBF\cite{Dehgan8} \cite{Dehgan8} & CS-RBF \cite{Dehgan8} & Optimal explicit \cite{Dehgan8}\\
 & &error & error &  error &  error &  error  \\
\hline

0.1  & 0.67249851 & $1.3 \times 10^{-5}$ & $1.5 \times 10^{-5}$ & $1.9 \times 10^{-5}$ &$1.6 \times 10^{-4}$ 
 & $5.2 \times 10^{-5}$\\
0.2  & 0.57206140 & $1.4 \times 10^{-5}$ & $1.5 \times 10^{-5}$ & $2.2 \times 10^{-5}$ & $1.9 \times 10^{-4}$ &$5.1 \times 10^{-5}$\\
0.3 & 0.41562694 & $1.0 \times 10^{-6}$  & $2.2 \times 10^{-6}$ & $2.3 \times 10^{-6}$ & $1.8 \times 10^{-5}$ &$5.1 \times 10^{-5}$\\ 
0.4   & 0.21850801 & $5.5 \times 10^{-6}$ & $2.2 \times 10^{-6}$ & $6.7 \times 10^{-7}$ & $8.0 \times 10^{-6}$ & $5.3 \times 10^{-5}$\\       
0.5  & 0.00000000 & $1.3 \times 10^{-7}$  & $1.3 \times 10^{-8}$ &$2.7 \times 10^{-7}$   & $8.7 \times 10^{-10}$ & $5.0 \times 10^{-5}$\\
0.6  & -0.21850801 & $5.2 \times 10^{-6}$ &  $1.3 \times 10^{-6}$ &$1.2 \times 10^{-6}$  & $8.0 \times 10^{-6}$ & $5.2 \times 10^{-5}$\\
0.7 & -0.41562694 & $1.0 \times 10^{-6}$ &  $2.2 \times 10^{-6}$ & $2.7 \times 10^{-6}$ & $1.8 \times 10^{-5}$ & $5.4 \times 10^{-5}$\\
0.8 & -0.57206140 & $1.4 \times 10^{-5}$ &  $1.5 \times 10^{-5}$ & $2.2 \times 10^{-5}$ &$1.9 \times 10^{-4}$  & $5.3 \times 10^{-5}$\\
0.9 & -0.67249851 & $1.3 \times 10^{-5}$ &  $1.5 \times 10^{-5}$ & $1.7 \times 10^{-5}$  &$1.6 \times 10^{-4}$  &$5.5 \times 10^{-5}$\\
1.0 & -0.70710678& $2.4 \times 10^{-7}$ &  $4.3 \times 10^{-9}$ & $2.1 \times 10^{-6}$ & $2.8 \times 10^{-9}$ &$5.4 \times 10^{-5}$\\ 
\hline

\end{tabular}
\caption{Point wise absolute error at $T = 0.25, \Delta t = 10^{-4}$ and $J$ = 6. }
\label{Table:3}
\end{table}

\begin{table}[H]
\centering

\begin{tabular}{|c|c|c|c|c|} 
\hline
$\Delta t$ & HHWCM &  TPS-RBF \cite{Dehgan8} & MQ-RBF \cite{Dehgan8} & CS-RBF \cite{Dehgan8}\\
\hline
$10^{-3}$  & $1.9 \times 10^{-5}$ & $6.8 \times 10^{-5}$ & $7.3 \times 10^{-5}$ & $2.3 \times 10^{-4}$\\ 
\hline
$10^{-4}$  & $1.4 \times 10^{-5}$ & $1.5 \times 10^{-5}$ & $2.2 \times 10^{-5}$ & $1.8 \times 10^{-4}$\\ 
\hline
\end{tabular}
\caption{Comparison of maximum error using various numerical methods at different time steps. }
\label{Table:4}
\end{table}

\section{Conclusion}
In this paper, we have developed a hybrid Haar wavelet collocation method for the numerical solution of  nonlocal hyperbolic partial differential equations. Instead of reformulating the original problem into periodic problem, we dealt with the integral boundary condition directly using the given data which is more accurate. For the spatial discretization, Haar wavelets are used whereas second order finite difference is used for temporal discretization. Stability analysis based on eigenvalue properties is carried out. We have derived error estimate for the proposed method. Finally, numerical results are presented and it is shown that our method is better than few existing method. This method can easily be generalized to higher dimensional problems.

\end{document}